\documentclass[11pt]{article}

\usepackage{amsfonts,graphicx,amsmath,amssymb,amsthm,geometry,bbm,hyperref,nicefrac,enumerate,comment,fontawesome,mathtools,cases,algorithmic,algorithm2e,stmaryrd,pict2e,picture,mathrsfs}
\usepackage[cal=boondoxo]{mathalpha} 

\usepackage{color, colortbl}
\definecolor{Gray}{gray}{0.9}

\newcommand{\R}{\mathbb{R}}
\newcommand{\N}{\mathbb{N}}
\newcommand{\E}{\mathbb{E}}

\renewcommand{\P}{\mathbb{P}}

\newcommand{\F}{\mathcal{F}}

\newcommand{\e}{\text{e}}
\renewcommand{\d}{\text{d}}

\usepackage{array}
\newcolumntype{P}[1]{>{\centering\arraybackslash}p{#1}}
\newcommand{\triple}{{\vert\kern-0.25ex\vert\kern-0.25ex\vert}}
\usepackage{mathtools}

\usepackage{geometry}
\usepackage{xcolor}

\newtheorem{remark}{Remark}[section]

\newtheorem{theorem}{Theorem}[section]

\newtheorem{example}{Example}[section]

\theoremstyle{definition}

\newcommand{\footremember}[2]{%
    \footnote{#2}
    \newcounter{#1}
    \setcounter{#1}{\value{footnote}}%
}
\newcommand{\footrecall}[1]{%
    \footnotemark[\value{#1}]%
} 

\begin{document}

\title{The deep multi-FBSDE method: a robust deep learning method for coupled FBSDEs}
\author{%
  Kristoffer Andersson\footremember{alley}{Mathematical Institute, Utrecht University, Utrecht, the Netherlands}\footremember{alleys3}{Department of Economics, University of Verona.
  Email: \href{mailto:kristofferherbert.andersson@univr.it}{kristofferherbert.andersson@univr.it}} 
  \and Adam Andersson \footremember{trailer}{Research Group of Computational Mathematics, Chalmers University of Technology \& University of Gothenburg}\footremember{alleys2}{Saab AB, Gothenburg, Sweden} \and Cornelis W. Oosterlee\footrecall{alley} 
 }


\maketitle
\begin{abstract}
We introduce the deep multi-FBSDE method for robust approximation of coupled forward–\\backward stochastic differential equations (FBSDEs), focusing on cases where the deep BSDE method of Han, Jentzen, and E (2018) fails to converge. To overcome the convergence issues, we consider a family of FBSDEs that are equivalent to the original problem in the sense that they satisfy the same associated partial differential equation (PDE) and initial value. Our algorithm proceeds in two phases: first, we approximate the initial condition jointly for a small number of FBSDEs from the FBSDE family, and second, we approximate the original FBSDE using the initial condition approximated in the first phase. Numerical experiments show that our method converges even when the standard deep FBSDE method does not.

\end{abstract}

\section{Introduction}
Solving forward-backward stochastic differential equations (FBSDEs) using deep learning has emerged as a powerful technique to overcome the curse of dimensionality inherent in classical numerical methods. Han, Jentzen, and E's seminal work~\cite{han2018solving} introduced the \emph{deep BSDE method}, employing neural networks to approximate solutions of high-dimensional parabolic partial differential equations (PDEs) by reformulating them as FBSDEs. This breakthrough demonstrated remarkable performance on problems previously considered intractable, showcasing the ability of neural networks to efficiently handle complex nonlinearities at scales unattainable by traditional approaches. Following this success, several deep learning-based strategies have emerged, notably \cite{beck2021solving,beck2019machine,beck2021deep,fujii2019asymptotic,raissi2018forward,ji2020three,ji2021control,wang2022deep,henry2017deep,andersson2022deep}, with convergence analyses provided in, e.g., \cite{han2020convergence,hutzenthaler2020proof,berner2020analysis,elbrachter2021dnn,Grohs2018APT,jentzen2018proof,reisinger2024posteriori,jiang2021convergence,andersson2025multi}. Concurrently, a separate branch known as backward-type methods, closer in spirit to classical dynamic programming algorithms, has developed, see e.g., \cite{hure2020deep,chan2019machine,fang2009novel,Balint,kapllani2024backward,hure2020deep,germain2022approximation}. For a comprehensive overview of machine learning algorithms for PDE approximation, we refer to \cite{beck2020overview}.

Despite these advancements, forward-type deep BSDE methods, including the original method in \cite{han2018solving}, still face significant convergence challenges when addressing certain classes of FBSDEs, such as those encountered in stochastic control problems~\cite{andersson2023convergence} and beyond. Specifically, these forward-type approaches often fail to converge or become trapped in suboptimal solutions due to ill-conditioned optimization landscapes~\cite{hure2020deep,andersson2023convergence,bussell2023deep}. Although backward-type methods may not suffer from precisely the same convergence issues, they encounter different challenges, especially implementation challenges when dealing with coupled FBSDEs (see \cite{chessari2021numerical,bender2008time,huang2025convergence2} for further discussions of implementation challenges for backward type methods when dealing with coupled FBSDEs) and accumulated errors in general (see e.g., \cite{Balint}). For FBSDEs related to stochastic control problems through dynamic programming, a robust deep FBSDE method was proposed in our paper \cite{andersson2023convergence}.  The scope of the present paper is the introduction of another robust modification of the deep FBSDE method, the \emph{deep multi-FBSDE method}. It is applicable to a broader class of FBSDEs and their PDEs, not strictly defined due to the experimental nature of our paper, but not restricted to only control problems. By robustness, we mean that the method can be optimized for the equation of interest, and not the other way around, which is sometimes seen in the deep learning-based scientific computing literature. Of course, this is an inprecise definition and no method allows convergence for any equation. Still, our method has been working on the vast majority of equations we have tested and we have deliberately chosen challenging equations. We refer to the original deep BSDE method, applied directly to coupled FBSDE without modification, as \emph{the deep FBSDE method} \cite{andersson2023convergence, ji2020three}.

Our proposed approach alleviates optimization difficulties by utilizing a family of FBSDEs that share the same initial condition as well as the same functional forms as the original FBSDE. In other words, all FBSDEs in this family are solved by the same PDE. We propose a two-phase optimization procedure. In the first phase, we approximate the joint initial conditions for all FBSDEs. In the second phase, we apply a reduced version of the deep FBSDE method, where we optimize only over the control processes while using the initial condition obtained from the first phase. This leads to a numerical scheme which we experimentally show has consistent and accurate convergence, even for strongly coupled, nonlinear FBSDEs. Our numerical experiments demonstrate stable and accurate solutions for problems for which the classical deep FBSDE methods fail. 

In summary, the primary contributions of this paper are as follows.
\begin{itemize}
\item A robust deep FBSDE method that we empirically demonstrate solves the convergence issues faced by existing deep learning-based methods for strongly coupled FBSDEs.
\item A general purpose training strategy that is demonstrated numerically to be applicable to general FBSDE systems, without restrictive assumptions such as the existence of an equivalent stochastic control problem, as in \cite{andersson2023convergence}.
\item Comprehensive numerical demonstrations illustrating the effectiveness and robustness of our method, especially in challenging scenarios involving stochastic control and nonlinear PDEs, thus underscoring its broad applicability and efficiency.
\end{itemize}

The rest of this paper is organized as follows: Section 2 presents the FBSDE problem and its solution by the non-linear Feynman--Kac formula through the PDE. Section 3 introduces our proposed robust deep learning method, detailing its formulation and implementation. Numerical results demonstrating its advantages are presented in Section 4, and we end in Section 5 with conclusions and future directions.

\section{Forward backward stochastic differential equations}
\label{sec2}

Throughout this paper, we let $T\in(0,\infty)$, $d,k\in\N$, $x_0\in\R^d$, $(W_t)_{t\in[0,T]}$ be a $k$-dimensional standard Brownian motion on a filtered probability space $(\Omega,\F,(\F_t)_{t\in[0,T]},\P)$ and the problem coefficients be denoted by $b\colon [0,T]\times\R^d\times\R\times\R^k\to\R^d$, $\sigma\colon[0,T]\times\R^d\to\R^{d\times k}$, $g\colon\R^d\to\R$ and $f\colon[0,T]\times\R^d\times\R\times\R^k\to\R$. The general form of the FBSDEs considered is given by the coupled system of equations
\begin{equation}
    \label{eq:FBSDE}\begin{dcases}
        X_t=x_0 + \int_0^tb(s,X_s,Y_s,Z_s)\d s +\int_0^t\sigma(s,X_s)\d W_s,\\
        Y_t = g(X_T) + \int_t^Tf(s,X_s,Y_s,Z_s)\d s - \int_t^T \langle
        Z_s,\sigma(s,X_s)\d W_s \rangle,\quad t\in[0,T].
    \end{dcases}
\end{equation}
The associate PDE is then a semilinear parabolic PDE on the form 
\begin{equation}\label{eq:PDE}\begin{dcases}
    \frac{\partial v}{\partial t}(t,x) + \mathcal{L}v(t,x) + f\big(t,x,v(t,x),\mathrm{D}_x v(t,x)\sigma(t,x)\big)=0, & (t,x)\in [0,T)\times \R^d,\\
    v(T,x)=g(x), & x\in \R^d.
    \end{dcases}
\end{equation}
Here $\mathcal{L}$ is the infinitesimal generator, acting on a twice continuously differentiable functions $\rho$ by:
\begin{align*}
\mathcal{L}\rho(x)=&b\big(t,x,\rho(t,x),\mathrm{D}_x\rho(t,x)\sigma(t,x)\big)^\top\mathrm{D}_x\rho(x)+\frac{1}{2}\text{Tr}\big(a(t,x)\mathrm{D}_x^2\rho(x)\big)\\
=&\sum_{i=1}^db_i\big(t,x,\rho(t,x),\mathrm{D}_x\rho(t,x)\sigma(t,x)\big)\partial_i\rho(x) + \frac{1}{2}\sum_{i,j=1}^da_{ij}(t,x)\partial_{i,j}\rho(x),
\end{align*}
where $a=\sigma\sigma^\top$ and $\text{Tr}(\cdot)$ denotes the trace of a matrix. The connection between \eqref{eq:FBSDE} and \eqref{eq:PDE} is established through the nonlinear Feynman--Kac formula, see \textit{e.g.,} \cite{zhang2017backward}, which states that 
\begin{equation*}
    Y_t=v(t,X_t),\  \text{and} \  Z_t=\mathrm{D}_xv(t,X_t).
\end{equation*}
Note that in the above, we have implicitly assumed that $v$ is differentiable with respect to the spatial variable. Although this holds true when \eqref{eq:PDE} possesses a classical solution, we can often only expect a solution in the viscosity sense. In this situation, derivatives of $v$ with respect to the space variable should be interpreted as set-valued sub-derivatives. Nevertheless, from our perspective, which is to find approximations of \eqref{eq:FBSDE} and \eqref{eq:PDE}, we can overlook this obstacle and view $Z$ as defined above.

\section{The deep FBSDE and deep multi-FBSDE methods}
In this section, we first introduce the classical deep FBSDE method and briefly discuss its convergence problems for coupled FBSDEs. Most importantly, we introduce our novel modification of the scheme. We use a Linear Quadratic (LQ) problem to demonstrate both the convergence problems of the deep FBSDE method and the convergence of our method. The details of the LQ-problem can be found in Section~\ref{lqgcp_num} together complementary experimental results.

\subsection{The deep FBSDE method}\label{sec:deepBSDE}
The deep FBSDE method relies on a reformulation of the FBSDE \eqref{eq:FBSDE} into two forward SDEs, one with an a priori unknown initial value.  It relies moreover on the Markov property of the FBSDE, which guarantees that $Z_t=\zeta^*(t,X_t)$, for some function $\zeta^*\colon[0,T]\times\R^d\to\R^k$, that we refer to as a Markov map. Optimization is done with respect to such functions and the initial values $y_0$ with the objective to satisfy the terminal condition. More precisely, the FBSDE \eqref{eq:FBSDE} is reformulated into the following variational problem
\begin{equation}\label{var0_FBSDE}\begin{dcases}
\underset{y_0,\zeta}{\mathrm{minimize}}\ 
 \E|Y_T^{y_0,\zeta}-g(X_T^{y_0,\zeta})|^2,\quad \text{where}\\
X_t^{y_0,\zeta}=x_0+\int_0^tb(s,X_s^{y_0,\zeta},Y_s^{y_0,\zeta},Z^{y_0,\zeta}_s)\d s + \int_0^t\sigma(s,X_s^{y_0,\zeta})\d W_s,\\
Y_t^{y_0,\zeta}= y_0-\int_0^tf(s,X_s^{y_0,\zeta},Y_s^{y_0,\zeta},Z^{y_0,\zeta}_s)\d s +\int_0^t \langle
Z^{y_0,\zeta}_s,\sigma(s,X_s^{y_0,\zeta})\d W_s\rangle,\\
Z^{y_0,\zeta}_t=\zeta(t,X_t^{y_0,\zeta}),\quad t\in[0,T].
    \end{dcases}
\end{equation}
Here $y_0$ and $\zeta$ are sought in appropriate spaces. It is clear that setting $y_0=Y_0$ and $\zeta=\zeta^*$ yields a minimizer of the objective function which attains zero and hence, solves \eqref{var0_FBSDE}. Moreover, assuming well-posedness and sufficient regularity of \eqref{eq:FBSDE}, existence and uniqueness of the minimizer are guaranteed. The connection between \eqref{eq:FBSDE} and \eqref{eq:PDE} yields $y_0^*=Y_0=v(0,x_0)$ and $\zeta^*=Z=\mathrm{D}_x v$.

We obtain a semi-discretized version of \eqref{var0_FBSDE} by approximating the It\^o integrals using the Euler--Maruyama scheme. We assume an equidistant temporal grid with $N\in\mathbb{N}$ time steps given by $t_0<t_1<\ldots<t_N$ with $t_0=0$, $t_N=T$ and step size $h=t_{n+1}-t_n$, and denote the \textit{i.i.d.} Wiener increments by $\Delta W_n=W_{t_{n+1}}-W_{t_n}\sim\mathcal{N}(0,h)$. The semi-discrete version of \eqref{var0_FBSDE} is given by:
\begin{equation}\label{var0_discrete_FBSDE}\begin{dcases}
\underset{y_0,\zeta=(\zeta_0,\ldots,\zeta_{N-1})}{\mathrm{minimize}}\ 
 \E\big|Y_N^{h,y_0,\zeta}-g(X_N^{h,y_0,\zeta})\big|^2,\quad \text{where}\\
X_{n+1}^{h,y_0,\zeta}=X_n^{h,y_0,\zeta}+b\big(t_n,X_n^{h,y_0,\zeta},Y_n^{h,y_0,\zeta},Z_n^{h,y_0,\zeta}\big)h + \sigma(t_n,X_n^{h,y_0,\zeta})\Delta W_n,\\
Y_{n+1}^{h,y_0,\zeta}= Y_n^{h,y_0,\zeta}-f\big(t_n,X_n^{h,y_0,\zeta},Y_n^{h,y_0,\zeta},Z_n^{h,y_0,\zeta}\big)h +\big\langle Z_n^{h,y_0,\zeta},\sigma\big(t_n,X_n^{h,y_0,\zeta}\big)\Delta W_n\big\rangle,\\
Z_n^{h,y_0,\zeta}=\zeta_n(X_n^{h,y_0,\zeta}),\  Y_0^{h,y_0,\zeta}=y_0,\ X_0^{h,y_0,\zeta}=x_0, \quad n\in\{0,1,\ldots,N-1\}.
    \end{dcases}
\end{equation}
In the deep FBSDE method, we search for $y_0\in\mathbb{R}$, a simple trainable parameter, and for $\zeta=(\zeta_n)_{n=0}^{N-1}$ in a parametric space defined by the architecture of a specified neural network.

In \cite{han2020convergence} an a posteriori error analysis of \eqref{var0_discrete_FBSDE} is provided under certain, rather restrictive assumptions. The full simulation error for the FBSDE can be bounded (up to a constant) by the time step $h$ and an a posteriori term given by the mean-squared error of the discrete terminal condition. Concretely, there exists a constant $C$, independent of $h$, such that for sufficiently small $h$
\begin{equation}\label{bound}
        \sup_{t\in[0,T]}\E|X_t-\hat{X}_t^{h,y_0,\zeta}|^2 + \sup_{t\in[0,T]}\E|Y_t-\hat{Y}_t^{h,y_0,\zeta}|^2 + \E\bigg[\int_0^T|Z_t-\hat{Z}_t^{h,y_0,\zeta}|^2\mathrm{d} t\bigg]\leq C\big(h + \E|\hat{Y}_N^{h,y_0,\zeta} - g(X_N^{h,y_0,\zeta})|^2),
    \end{equation}
where for $n\in\{0,1,\ldots,N-1\}$, and $t\in[t_n,t_{n+1})$ we have defined
\begin{equation*}
    \hat{X}_t^{h,y_0,\zeta} \;=\; X_n^{h,y_0,\zeta},
    \quad
    \hat{Y}_t^{h,y_0,\zeta} \;=\; Y_n^{h,y_0,\zeta},
    \quad
    \hat{Z}_t^{h,y_0,\zeta} \;=\; Z_n^{h,y_0,\zeta}.
\end{equation*}
Under conditions that guarantee the same or similar a posteriori bound, such as those in \cite{han2020convergence}, the function $\E|\hat{Y}_N^{h,y_0,\zeta} - g(X_N^{h,y_0,\zeta})|^2$ is a highly suitable loss function and guarantees an approximation with an error bounded by the loss after training and the time step. When the assumptions are not met, the situation is less fortunate. Counterexamples of FBSDEs from stochastic control are presented in \cite{andersson2023convergence} for which it is possible to make the right-hand side of \eqref{bound} arbitrarily small, but with the left-hand side remaining large and discrete solutions converging to severely wrong processes. In this paper, we provide four more counterexamples. Thus, while the a posteriori bound in \cite{han2020convergence} is a powerful tool under its specific conditions, extrapolating it beyond these conditions is a questionable practice.

\begin{figure}[htp]
\centering
\begin{tabular}{cc}
          \includegraphics[width=80mm]{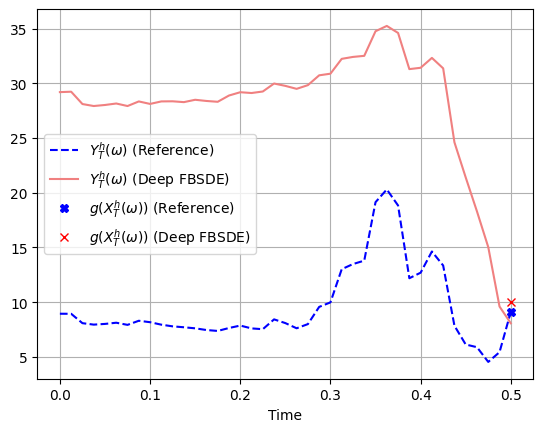}&    
          \includegraphics[width=80mm]{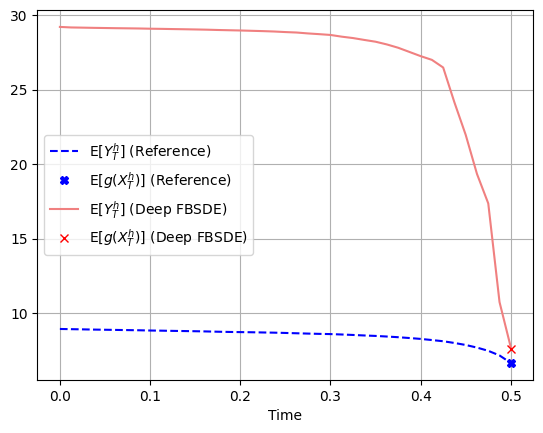}
          \end{tabular}
\caption{The deep FBSDE method failing to approximate the true reference solution. \textbf{Left:} One representative path. \textbf{Right:} A sample mean of size $2^{12}$.} \label{YEY}
\end{figure}

For instructive purposes and since we address the convergence problem documented in \cite{andersson2023convergence}, with a novel method, we briefly repeat the problem with a numerical example. In order to demonstrate that the optimization landscape for the deep FBSDE method is problematic for strongly coupled FBSDEs, we use the mean-squared error (MSE)
\begin{equation}\label{var_empirical_investigation}
\mathrm{MSE}(y_0)
\coloneqq
\underset{\zeta=(\zeta_0,\dots,\zeta_{N-1})}{\inf}\
\E
\big[
\big|
Y_N^{h,y_0,\zeta}-g(X_N^{h,y_0,\zeta})
\big|^2
\Big].
\end{equation}
This is the MSE with $y_0$ fixed and only $\zeta$ being optimized. If joint optimization of $y_0$ and $\zeta$, as in the deep FBSDE method, is supposed to yield a good approximation of the FBSDE, then by necessity the minimum of $y_0\mapsto\mathrm{MSE}(y_0)$ must be close to the true $Y_0$ of the FBSDE. For an LQ problem, detailed in Section~\ref{lqgcp_num} below, Figure~\ref{fig:MSE} shows that the optimization landscape $y_0\mapsto\mathrm{MSE}(y_0)$ instead has a minimum $y_0^*$ very far from the true $Y_0$. Figure~\ref{YEY} shows the exact $Y$ and the output of the deep FBSDE method, pathwise and in mean. We note that although the initial condition is significantly incorrect, the terminal condition is well approximated, i.e., $g(X_N^{h,y_0,\zeta})-Y_N^{h,y_0,\zeta}$ is small both on the pathwise level and in the mean-squared sense. It should be noted that the numerical algorithm has converged, as the loss function is no longer decreasing and $y_0$ remains constant (not shown here but Figure~\ref{fig:loss_CHT} qualitatively shows the same behavior for a different problem). Moreover, the algorithm is consistent, with repeated runs yielding similar values for the loss function and $y_0$. The reader might wonder if a smaller time step $h$ in the approximation would improve approximation of $Y_0$. Figure~\ref{fig:L2_Y0} indicates a negative answer to this valid thought by showing non-convergence of our implementation of the deep FBSDE method (after thorough hyperparameter optimization) (the case $K=1$ in that figure).

\begin{figure}[htp]
\centering
\begin{tabular}{c}
          \includegraphics[width=120mm]{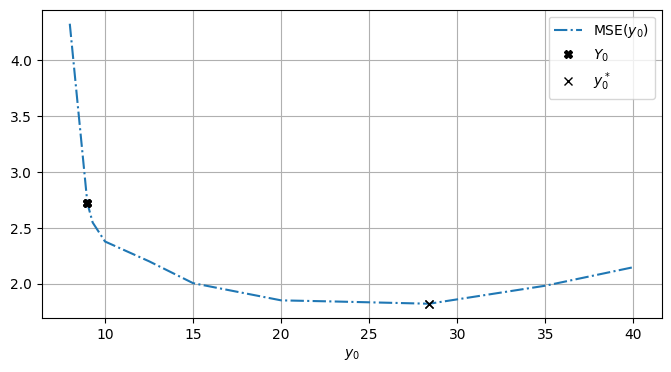}
          \end{tabular}
\caption{MSE plotted for different values of $y_0$ for the LQ-problem of Section~\ref{lqgcp_num}.}\label{fig:MSE}
\end{figure} 

The behavior depicted in Figures~\ref{fig:MSE}--\ref{YEY} is not confined to a narrow class of specific problems or special cases. For example, it is not limited to strongly coupled FBSDEs, and one can easily construct low-dimensional (e.g., $d=2$) examples of decoupled BSDEs with smooth solutions to the associated PDE, wherein the deep FBSDE method fails. In particular, quadratic BSDEs appear  susceptible to this phenomenon. A more detailed, although mostly qualitative and empirical, discussion of the reasons for this lack of convergence is provided in \cite{andersson2023convergence}.

\subsection{The deep multi-FBSDE method}

Our deep multi-FBSDE method is based on a slightly modified variational formulation. To arrive at it, we notice that for sufficiently regular $\psi\colon[0,T]\times\R^d\times\R\times\R^\ell\to\R^d$, it holds by the It\^o formula, for $Y_t^\psi=v(t,X_t^\psi)$ and $Z_t^\psi=\mathrm{D}_xv(t,X_t^\psi)$ that
\begin{equation}\label{eq:FBSDE_drift}\begin{dcases}
  X_t^\psi
  =
  x_0+\int_0^t \big(b(s,X_s^\psi,Y_s^\psi, Z_s^\psi)-\psi(s,X_s^\psi,Y_s^\psi,Z_s^\psi)\big)\d s + \int_0^t\sigma(s,X_s^\psi)\d W_s,\\
  Y_t^\psi
  =g(X_T^\psi) 
  + 
  \int_t^T 
    \big(f(s,X_s^\psi,Y_s^\psi,Z_s^\psi) +\big\langle \psi(s,X_s^\psi,Y_s^\psi,Z_s^\psi),Z_s^\psi\big\rangle\big)
  \d s \\ \quad\quad\ 
  - \int_t^T  \big\langle Z_s^\psi ,\sigma(s,X_s^\psi) \d W_s\big\rangle,\quad t\in[0,T].
    \end{dcases}
\end{equation}
This defines a family of FBSDEs, indexed by $\psi$, which are all equivalent to the PDE~\eqref{eq:PDE}. In particular, for all $\psi$ we have $Y_0^\psi=v(0,X_0^\psi)=v(0,x_0)$ and thus the initial value is shared for all $\psi$, which we utilize. Introducing the shorthand notation
\begin{align*}
  b^{\psi}(t,x,y,z)
  &:=
  b(t,x,y,z) - \psi(t,x,y,z)\\
  f^{\psi}(t,x,y,z)
  &:=
  f(t,x,y,z)+\langle z, \psi(t,x,y,z)\rangle,
\end{align*}
we have for any $\psi$ the variational formulation 
\begin{equation}\label{var_Psi_FBSDE}\begin{dcases}
\underset{y_0,\zeta}{\mathrm{minimize}}\ 
 \E\big|Y_T^{y_0,\zeta,\psi}-g(X_T^{y_0,\zeta,\psi})\big|^2,\quad \text{where} \\
X_t^{y_0,\zeta,\psi}=x_0+\int_0^t b^\psi(s,X_s^{y_0,\zeta,\psi},Y_s^{y_0,\zeta,\psi},Z^{y_0,\zeta,\psi}_s)\d s  + \int_0^t\sigma(s,X_s^{y_0,\zeta,\psi})\d W_s,\\
Y_t^{y_0,\zeta,\psi}= y_0-\int_0^t f^\psi(s,X_s^{y_0,\zeta,\psi},Y_s^{y_0,\zeta,\psi},Z^{y_0,\zeta,\psi}_s)\d s +\int_0^t \big\langle Z^{y_0,\zeta,\psi}_s,\sigma(s,X_s^{y_0,\zeta,\psi})\d W_s \big\rangle,\\
Z^{y_0,\zeta,\psi}_t=\zeta(t,X_t^{y_0,\zeta,\psi}),\quad t\in[0,T].
    \end{dcases}
\end{equation}
Since the optimal $y_0^*$ and $\zeta^*$ are optimal for all $\psi$, we trivially have, for $K\in\N$, $\psi_1,\dots,\psi_K\colon[0,T]\times\R^d\times\R\times\R^\ell\to\R^d$, the equivalent variational formulation,
\begin{equation}\label{var_Psi_FBSDE2}\begin{dcases}
\underset{y_0,\zeta}{\mathrm{minimize}}\ 
 \sum_{k=1}^K\E\big|Y_T^{y_0,\zeta,\psi_k}-g(X_T^{y_0,\zeta,\psi_k})\big|^2,\quad \text{where for }\psi\in\{\psi_1,\dots,\psi_K\} \\
X_t^{y_0,\zeta,\psi}=x_0+\int_0^t b^\psi(s,X_s^{y_0,\zeta,\psi},Y_s^{y_0,\zeta,\psi},Z^{y_0,\zeta,\psi}_s)\d s  + \int_0^t\sigma(s,X_s^{y_0,\zeta,\psi})\d W_s,\\
Y_t^{y_0,\zeta,\psi}= y_0-\int_0^t f^\psi(s,X_s^{y_0,\zeta,\psi},Y_s^{y_0,\zeta,\psi},Z^{y_0,\zeta,\psi}_s)\d s +\int_0^t \big\langle Z^{y_0,\zeta,\psi}_s,\sigma(s,X_s^{y_0,\zeta,\psi})\d W_s\big\rangle,\\
Z^{y_0,\zeta,\psi}_t=\zeta(t,X_t^{y_0,\zeta,\psi}),\quad t\in[0,T].
    \end{dcases}
\end{equation} 

The deep multi-FBSDE is divided into two phases. Phase I is based on approximating \eqref{var_Psi_FBSDE2} for suitable $K$ and $\psi_1,\dots,\psi_K$, with time discretization, neural network approximation of $\zeta$ and stochastic gradient descent. In the second phase, a better approximation of $\zeta$ is obtained by solving the classical deep FBSDE method with the exception that only $\zeta$ is optimized and $y_0$ is taken from Phase~I. For easy and low-dimensional problems, Phase~I gives an accurate approximation of $\zeta$ and Phase~II is not needed. For more challenging and high-dimensional problems, we have observed that Phase~II improves the approximation accuracy. We believe that this has to do with the fact that in Phase~I the Markov map $\zeta$ is optimized for good fit along typical trajectories of $K$ different SDEs, i.e., in a larger domain than in Phase~II, where it is only trained along typical trajectories of the original forward SDE, i.e., for a more limited task. 

We next state the scheme and similar to the deep FBSDE method let $\zeta$ represent $\zeta_0,\ldots,\zeta_{N-1}$.

\vspace{0.5cm}
\noindent\textbf{Phase I - Approximation of $Y_0$:}\newline
\begin{equation}\label{var_discrete_FBSDE_K}\begin{dcases}
\underset{y_0,\zeta_0,\ldots,\zeta_{N-1}}{\mathrm{minimize}}\ 
 \sum_{k=1}^K\E\big|Y_N^{h,y_0,\zeta,\psi_k}-g(X_N^{h,y_0,\zeta,\psi_k})\big|^2,\quad \text{where for } \psi\in\{\psi_1,\dots,\psi_K\}\\
X_{n+1}^{h,y_0,\zeta,\psi}=X_n^{h,y_0,\zeta,\psi}+b^{\psi}\big(t_n,X_n^{h,y_0,\zeta,\psi},Y_n^{h,y_0,\zeta,\psi},Z_n^{h,y_0,\zeta,\psi}\big)h + \sigma(t_n,X_n^{h,y_0,\zeta,\psi})\Delta W_n,\\
Y_{n+1}^{h,y_0,\zeta,\psi} =Y_n^{h,y_0,\zeta,\psi}-f^{\psi}\big(t_n,X_n^{h,y_0,\zeta,\psi},Y_n^{h,y_0,\zeta,\psi},Z_n^{h,y_0,\zeta,\psi}\big)h +\big\langle Z_n^{h,y_0,\zeta,\psi},\sigma\big(t_n,X_n^{h,y_0,\zeta,\psi}\big)\Delta W_n\big\rangle,\\
Z_n^{h,y_0,\zeta,\psi}=\zeta_n(X_n^{h,y_0,\zeta,\psi}\big),\  X_0^{h,y_0,\zeta,\psi}=x_0,\ Y_0^{h,y_0,\zeta,\psi}=y_0,\\
\text{for }n\in\{0,1,\ldots,N-1\}.
    \end{dcases}
\end{equation}

\vspace{0.5cm}
\noindent\textbf{Phase II - Approximation of the full FBSDE with known initial value:}\newline
\begin{equation}\label{var_discrete_FBSDE_fix_y}\begin{dcases}
\underset{\zeta_0,\ldots,\zeta_{N-1}}{\mathrm{minimize}}\ 
 \E\big|Y_N^{h,\zeta}-g(X_N^{h,\zeta})\big|^2,\quad \text{where}\\
X_{n+1}^{h,\zeta}=X_n^{h,\zeta}+b\big(t_n,X_n^h,Y_n^{h,\zeta},Z_n^h\big)h + \sigma(t_n,X_n^{h,\zeta})\Delta W_n,\\
Y_{n+1}^{h,\zeta}= Y_n^{h,\zeta}-f\big(t_n,X_n^{h,\zeta},Y_n^{h,\zeta},Z_n^{h,\zeta}\big)h +\big\langle Z_n^{h,\zeta},\sigma\big(t_n,X_n^{h,\zeta}\big)\Delta W_n\big\rangle,\\
Z_n^{h,\zeta}=\zeta_n\big(X_n^{h,\zeta}\big),\  Y_0^{h,\zeta}=y_0^*, \quad \text{for } n\in\{0,1,\ldots,N-1\}.
    \end{dcases}
\end{equation}

\begin{figure}[htp]
\centering
\begin{tabular}{c}
          \includegraphics[width=120mm]{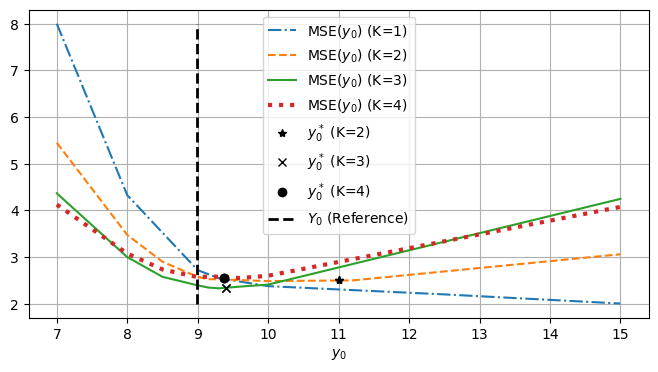}
          \end{tabular}
\caption{Four MSE curves plotted for $K=1,2,3,4$, with $\psi_1=0$, $(\psi_1,\psi_2)=(0,b_1)$, $(\psi_1,\psi_2,\psi_3)=(0,b_1,-b_1)$ and $(\psi_1,\psi_2,\psi_3,\psi_4)=(0,b_1,-b_1, -0.5b_1)$, where $b(x,y,z)=b_1(x)+b_2(z)$. In Phase II, $y_0^*$ for $K=1,2,3,4$ are the initial conditions approximated with Phase I. Note that for $K=1$, $y_0^*=28.44$ as can be visualized in Figure~\ref{fig:MSE}. For visualization purposes, we have scaled the curves. For $K=1,2,3,4$, the training times are approximately 400s, 500s, 600s, and 700s, respectively.}\label{fig:5MSE}
\end{figure}

\begin{figure}[htp]
\centering
\begin{tabular}{c}
          \includegraphics[width=120mm]{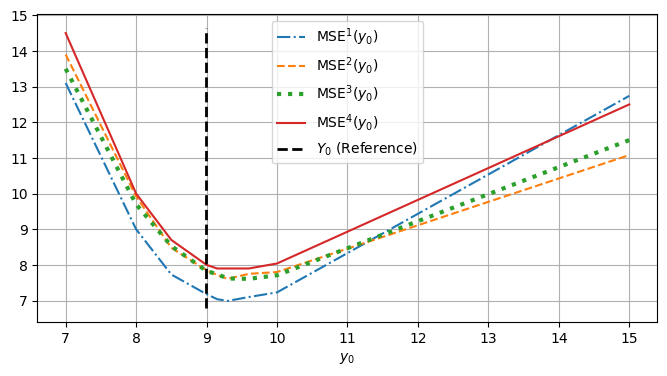}
          \end{tabular}
\caption{Four MSE curves plotted for $\psi_1=b_1\circ\xi$, $\psi_2=-b_1\circ\xi$, $\psi_3=0$, for $\xi\in\{x,-\sin(x),x\cos(x),1-\exp((x\wedge 1)\vee(-1))\}$, where $\xi$ acts coordinate-wise on vectors. The order coincides with that of the set.}\label{fig:4MSE}
\end{figure}

\begin{figure}[htp]
\centering
\begin{tabular}{c}
          \includegraphics[width=120mm]{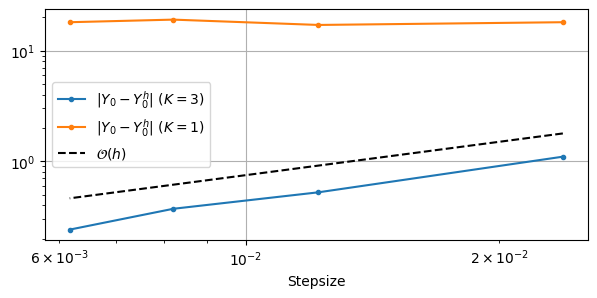}
          \end{tabular}
\caption{Empirical convergence plot for approximations of $Y_0$ using the deep multi-FBSDE method Phase~I. Here the case $K=1$ coincides with the original deep FBSDE method.}\label{fig:L2_Y0}
\end{figure}

We demonstrate the new numerical scheme on the LQG-problem, detailed in Section~\ref{lqgcp_num}, for which the deep FBSDE method was shown to fail in Section~\ref{sec:deepBSDE}. In Phase~I, we take $K=1,2,3,4$, with the choices $\psi_1=0$, $(\psi_1,\psi_2)=(0,b_1)$, $(\psi_1,\psi_2,\psi_3)=(0,b_1,-b_1)$ and $(\psi_1,\psi_2,\psi_3,\psi_4)=(0,b_1,-b_1,-0.5b_1)$, where the drift has the form $b(x,y,z)=b_1(x)+b_2(z)$, and plot the MSE curves in Figure~\ref{fig:5MSE}. For $K=2,3,4$ the optimization landscape is much more satisfactory than for $K=1$, which is the deep FBSDE method. In fact the minima are close to the real $Y_0=8.94$. To further indicate that the method is robust to the choice of $K\geq3$ (and possibly $K=2$) and $(\psi_1,\dots,\psi_K)$ we take $K=3$ and rather arbitrary choices $\psi_1=b_1\circ\xi$, $\psi_2=-b_1\circ\xi$, $\psi_3=0$, for $\xi\in\{x,-\sin(x),x\cos(x),1-\exp(\min(1,\max(-1,x)))\}$, where $\xi$ acts coordinate-wise on vectors. We obtain the four loss surfaces in Figure~\ref{fig:4MSE}, which, based on visual inspection, have the same minima. This robustness is a very important property because when the method is applied to a problem without a reference solution, one must be able to trust that the choice of $\psi_1,\psi_2,\psi_3$ does not need to be tuned based on comparison with a reference solution. Our results indicate a robustness in the choice of $K\geq2$ and $(\psi_1,\dots,\psi_K)$, but it needs to be investigated further, particularly theoretically. The deviation between $Y_0$ and the optimized $y_0^*$ is the discretization error due to the time stepping. This is seen empirically in the convergence plot in Figure~\ref{fig:L2_Y0} for $\xi(x)=x$, showing the experimental convergence order one, i.e., $|Y_0-y_0^*(h)|=\mathcal O(h)$. The convergence of $(X,Y,Z)$ as processes after Phase~II is shown in Section~\ref{lqgcp_num}.



\color{black}

\color{black}

\section{Numerical experiments}
We present four numerical experiments, two with an associated stochastic optimal control problem, and two without such a connection, all with analytic or semi-analytic reference solutions so that we can analyze the error. The first problem is a controlled Brownian motion in two dimensions, the second problem is a six-dimensional LQ-problem, and the last two examples are FBSDEs related to two-dimensional advection diffusion reaction PDEs, one without a reaction term and one with a linear one. We demonstrate the performance of both the deep FBSDE and the deep multi-FBSDE methods.

In the literature, it is most common to demonstrate deep BSDE-type methods on very high-dimensional problems. However, these problems are most often very symmetric, with solutions being permutation-invariant in its space variables. It is highly questionable how important these examples are for applications and also if they fairly represent more general problems of the same dimensions. We have instead chosen difficult low-dimensional problems for which the classical deep FBSDE methods fail, rather than high-dimensional problems.

In the implementation of the deep FBSDE and both Phase~I--II of the deep multi-FBSDE methods, we optimize $\zeta_0,\dots,\zeta_{N-1}$ in a function class defined by ReLU feed-forward neural networks with 3 hidden layers with 20 nodes in each layer. The initial value $y_0$ in the deep FBSDE method and in Phase~I of the deep multi-FBSDE method is represented by one single trainable parameter. Training data consist of $2^{20}$ independent realizations of the $d$-dimensional Wiener increment vectors $\Delta W_0,\ldots,\Delta W_{N-1}$, each element being \textit{i.i.d.}, and normally distributed with mean zero and variance $h$. The training data are reused in $10$ epochs. Training is initiated by randomly sampling parameters and proceeds with the classical Adam optimizer with batches of size $2^{12}$.

\subsection{FBSDEs stemming from stochastic optimal control problems}
In this section, we look into two main types of stochastic control problems, controlled Brownian motions and linear quadratic Gaussian control problems. Before discussing these specific topics, we define the general form of the stochastic control problem:
\begin{equation}\label{control_problem}
\begin{dcases}
\underset{u}{\mathrm{minimize }}\;\;
J(u;0,x),\quad \text{where}\\
\d X_t = b(t,X_t,u_t) \d t + \sigma(t,X_t) \d W_t, \\
J(u;t,x) = \E \bigg[\int_t^T f(s,X_s,u_s) \d s + g(X_T) \,\big|\, X_t=x\bigg], \quad t \in [0,T],\ x\in\R^d.
\end{dcases}
\end{equation}
The controls $u$ belong to a suitable subset of adapted stochastic processes, and we refrain from details  The value function associated with the control problem is formally defined as \begin{equation*}
    v(t,x)=\inf_{u}J(u;t,x).
    \end{equation*}
Under sufficient regularity assumptions on $b,\sigma,f,g$, a solution to the following Hamilton--Jacobi--Bellman (HJB) equation exists and coincides with the value function of the control problem
\begin{equation}\label{PDE_CH}\begin{dcases}
\frac{\partial v}{\partial t} + \frac{1}{2}\mathrm{Tr}\big(\sigma^\top \sigma \mathrm{Hess}_x v\big) +\inf_{u}\big[  \langle\mathrm{D}_x v,b\rangle +f\big] = 0, & (t,x) \in [0,T)\times \R^d,\\
v(T,x) = g(x), & x \in  \R^d.
\end{dcases}
\end{equation}
In this paper we tacitly assume such regularity. The resulting HJB-equation is a semilinear parabolic PDE of the type discussed in Section~\ref{sec2}; therefore, it admits an equivalent formulation as an FBSDE.

\subsubsection{Controlled Brownian motions with general terminal cost}
In this section, we set $k=d=\ell\in\mathbbm{N}$, $x_0 \in \R^d$,  and $\sigma,r \in \R^+$. For $t \in [0, T], x \in \R^d,$ and $u \in \R^\ell$, the drift, diffusion and cost fu
nctions are defined, as follows, 
\begin{equation}
b(t,x,u) = u, \quad\sigma(t,x)= \sigma I_d,\quad f(t,x,u) = \frac{r}{2}\|u\|^2,
\end{equation}
together with a flexibly chosen terminal cost function $g\colon\R^d\to\R$. Straightforward calculations yield the optimal control $
    u^*=-\frac{1}{r}\mathrm{D}_x v=-\frac{1}{r}Z$, from which we obtain the FBSDE
\begin{equation}\label{FBSDE_CHT}\begin{dcases}
X_t = x_0 - \int_0^t \frac{1}{r}Z_s \text{d} s + \sigma W_t,\\
Y_t = g(X_T) + \int_t^T \frac{1}{2r}\|Z_s\|^2 \text{d} s - \int_t^T \langle Z_s ,\sigma \text{d} W_s\rangle, \quad t\in[0,T].
\end{dcases}
\end{equation}
Using the Cole--Hopf transformation, the unique $Y$-component of the solution to the above FBSDE is
\begin{equation}\label{MC_PDE}
    Y_t=-r\sigma^2\log\bigg(\E\big[\text{e}^{-\frac{1}{r\sigma^2} g(X_t+\sigma\sqrt{T-t}\,\xi)}\big|X_t\big]\bigg),
    \quad
    t\in[0,T],
\end{equation}
where $\xi$ is a standard $d$-dimensional normally distributed random vector.

The FBSDE \eqref{FBSDE_CHT} can then be approximated using Monte Carlo sampling. It is important to note that, to approximate $(Y_t(\omega))_{t\in[0,T]}$, nested expectations are required along the trajectory of $(X_t(\omega))_{t\in[0,T]}$. Moreover, since $X$ is coupled with $Z$ in the drift component, a straightforward application of this strategy to approximate $Y$ at times other than $t=0$ (which is independent of the trajectory of $X$) is challenging. However, if our approximation of $(X,Y,Z)$ is accurate, it should follow that our approximation of $Y$ aligns with the Monte Carlo approximation of $Y$ using \eqref{MC_PDE}, along the simulated trajectories of $X$.
\newline\newline
\noindent \textbf{Problem specific settings:}\newline
For the deep multi-FBSDE method, we opt for $K=2$, setting $\psi_1(t,x,y,z)=0$ and $\psi_2(t,x,y,z)=\frac{r_1^2}{2r_2}z$. This configuration implies that $(X^{\psi_1},Y^{\psi_1},Z^{\psi_1})$ is the solution to \eqref{FBSDE_CHT}, whereas $(X^{\psi_2},Y^{\psi_2},Z^{\psi_2})$ satisfies the following decoupled equation:
\begin{equation}\label{FBSDE_CHT_decoupled}\begin{dcases}
X_t^{\psi_2} = x_0 + \sigma W_t,\\
Y_t^{\psi_2} = g(X_T^{\psi_2}) - \int_t^T \frac{1}{2r}\|Z_s^{\psi_2}\|^2 \text{d} s - \int_t^T \big\langle Z_s^{\psi_2}, \sigma \text{d} W_s\big\rangle, \quad t \in [0,T].
\end{dcases}
\end{equation}
For this relatively straightforward problem, we have observed that only the first phase of the deep multi-FBSDE method is necessary, so Phase II is skipped.

\begin{figure}[htp]
\centering
\begin{tabular}{ccc}
\includegraphics[width=80mm]{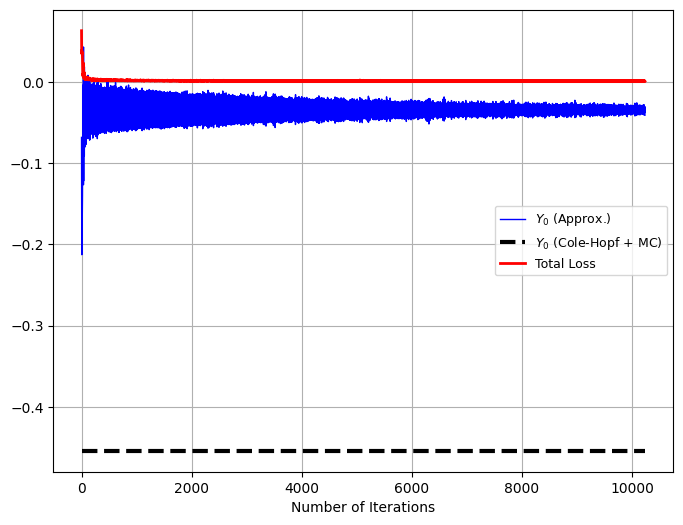}& \includegraphics[width=80mm]{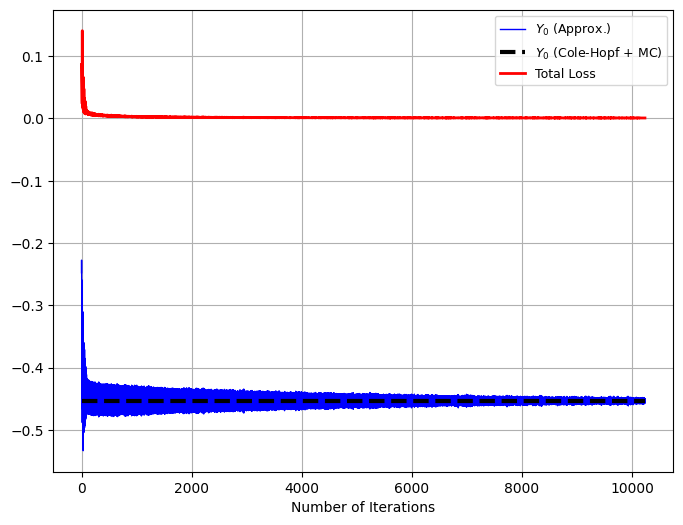}\
\end{tabular}
\caption{Loss and $Y_0$ during training for the deep FBSDE method (left) and the deep multi-FBSDE method (right). ''Number of iterations'' refers to gradient updates. The losses of the first two iterations were removed due to large values and for better visual representation.}\label{fig:loss_CHT}
\end{figure}

We set $T=0.5$, $d=\ell=k=2$, $r=1$, and $\sigma=0.25$. The initial state is set at $x_0=(-0.1, 0.1)^\top$. We employ a terminal cost function $g(x)= -|x_1-x_2|$, where $x=(x_1,x_2)^\top$. This configuration aims to control the two components of the Brownian motion such that they diverge as much as possible by time $t=T$, while minimizing the use of control force. The discretization of the problem is set with $N=40$. 
\newline\newline
\noindent \textbf{Results and discussion:}\newline
Figure~\ref{fig:loss_CHT} illustrates the loss and the approximation of $Y_0$ during training for both the deep FBSDE method and the deep multi-FBSDE method. Notably, while the loss approaches zero for both methods, only the deep multi-FBSDE method converges to the true $Y_0$.

\begin{figure}[htp]
\centering
\begin{tabular}{ccc}
\includegraphics[width=80mm]{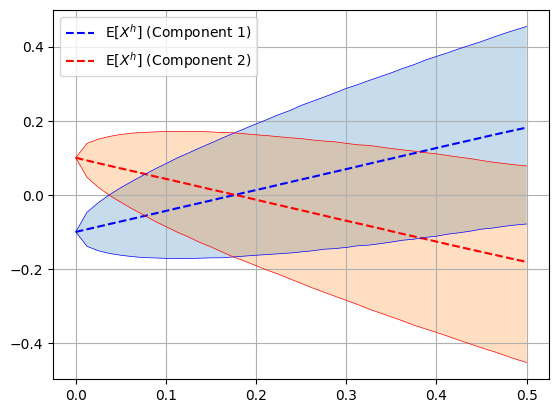}& \includegraphics[width=80mm]{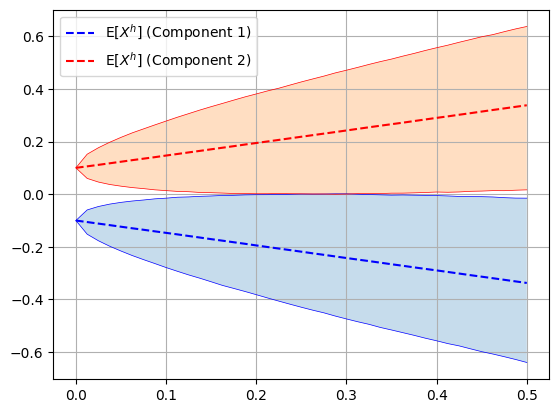}\\
\includegraphics[width=80mm]{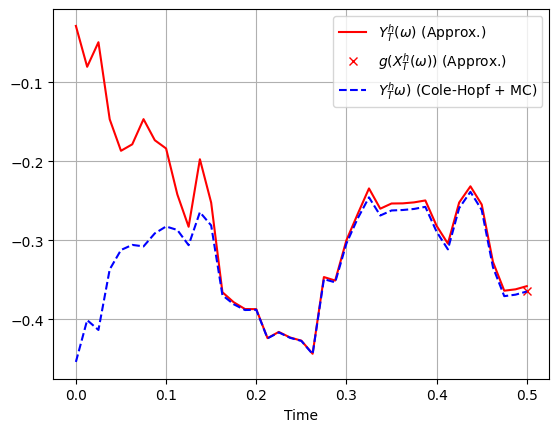}& \includegraphics[width=80mm]{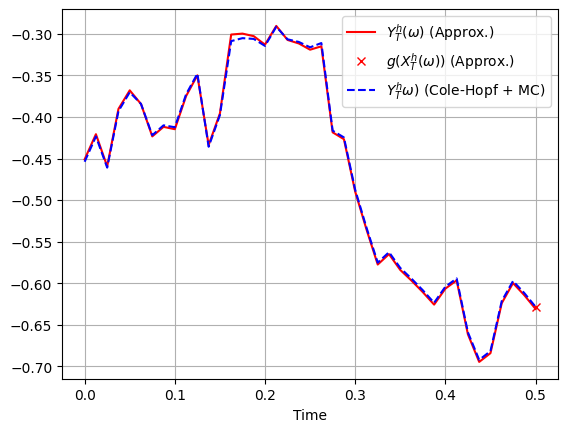}\
\end{tabular}
\caption{The upper plots display the empirical mean of $X$ as well as the 95th and 5th percentiles for the deep FBSDE method (left) and the deep multi-FBSDE method (right). The lower plots show a representative trajectory of the approximate $Y$ for the deep FBSDE method (left) and the deep multi-FBSDE method (right), compared with their semi-analytic counterparts.}\label{fig:FBSDE_CHT}
\end{figure}

The deep FBSDE method fails to converge to the true solution. Although we do not have theoretical understanding for the reason, to gain some partial insight, we illustrate in Figure~\ref{fig:FBSDE_CHT} the sample mean and a representative path of $X$ and $Y$, for the two numerical methods and the reference solution.
We see that the optimal strategy pushes $X^1$ to lower values and $X^2$ to higher values, which is intuitive since $X_0^1<X_0^2$ and the task is to maximize $|X_T^1-X_T^2|$ with low control cost. The deep FBSDE method, on the other hand, causes $X^1$ and $X^2$ to instead cross each other, and naturally uses a more expensive control for the task, reading the initially high values of the $Y$-process. After the crossing, the control seems to resemble the optimal control well, which can be read from the resulting $Y$-process being close to the reference. A similar MSE-analysis as in Figure~\ref{fig:MSE} (omitted) shows that the problem with the deep FBSDE solution is due to a bad global minimum.

\subsubsection{Linear Quadratic Gaussian control problems}\label{lqgcp_num}
Among all stochastic control problems, the linear quadratic Gaussian control problem is notably the most structured and extensively studied, see, for example, \cite{aastrom2012introduction}. For our purposes, it presents a closed-form analytic solution, allowing us to benchmark our numerical approximations against it. 

Consider the following settings: let $k=d\in\mathbbm{N}$, $\ell\in\mathbbm{N}$, $x_0\in\R^d$, $A, \sigma \in\R^{d\times d}$, $R_x, G\in\mathbbm{S}^d_+$, $R_u\in\mathbbm{S}^\ell_+$, and $B\in\R^{d\times \ell}$ with full rank, alongside $C\in\R^d$. For $t\in[0,T], x\in\R^d,$ and $ u\in\R^\ell$ let
\begin{equation*}
    b(t,x,u)=A(C-x) + Bu, \quad f(t,x,u) = \langle R_x x,x\rangle + \langle R_u u,u\rangle\quad g(x)=\langle Gx,x\rangle. 
\end{equation*}
The above, together with \eqref{control_problem}, defines the state equation and cost functional for a LQ control problem. The optimal feedback control, minimizing the Hamiltonian, $\inf_{u \in U}{\langle \text{D}_x V, Bu\rangle + \langle R_u u, u\rangle}$, is given by
\begin{equation}\label{opt_control}
u^*_t = -\frac{1}{2} R_u^{-1} B^\top \text{D}_x v(t, X_t).
\end{equation}
Recalling that $v$ is the solution to the associated HJB-equation, its solution can be expressed as
\begin{align}
v(t,x) = x^\top P(t)x + x^\top Q(t) + R(t),
\end{align}
where $(P, Q, R)$ are the solutions to the system of ordinary differential equations,
\begin{equation*}\begin{dcases}
\dot{P}(t) - A^\top P(t) - P(t)A - P(t)BR_u^{-1}B^\top P(t) + R_x = \mathbf{0}_{d\times d},\\
\dot{Q}(t) + 2P(t)AC - A^\top Q(t) - P(t)BR_u^{-1}B^\top Q(t) = \mathbf{0}_d,\\
\dot{R}(t) + \text{Tr}\big(\sigma\sigma^\top P(t)\big) + Q(t)^\top AC - \frac{1}{4}Q(t)^\top BR_u^{-1}B^\top Q(t) = 0, \quad t \in [0,T],\\
P(T) = G; \quad Q(T) = \mathbf{0}_d; \quad R(T) = 0.
\end{dcases}
\end{equation*}
We refer to the entire system colloquially as the Riccati equation, although the first equation is strictly a matrix Riccati equation. The gradient of $v$, $\text{D}_xv(t,x) = 2P(t)x + Q(t)$, leads to the associated FBSDE:
\begin{equation}
\begin{dcases}
X_t = x_0 + \int_0^t \big[A(C-X_s) - \frac{1}{2} BR_u^{-1}B^\top Z_s\big] \d s + \int_0^t \sigma \d W_s,\\
Y_t = \langle G X_T, X_T\rangle + \int_t^T \big(\langle R_x X_s, X_s\rangle + \frac{1}{4}\big\langle R_u^{-1} B^\top Z_s, B^\top Z_s\big\rangle\big) \d s - \int_t^T \langle Z_s, \sigma \d W_s\rangle, \quad t \in [0,T].
\end{dcases}\label{FBSDE_LQG}
\end{equation}
The solution to \eqref{FBSDE_LQG} is provided by
\begin{equation}\label{eq:YZ}
Y_t = X_t^\top P(t)X_t + X_t^\top Q(t) + R(t); \quad Z_t = 2P(t)X_t + Q(t).
\end{equation}
In our experiments, we use the Euler approximation of the Riccati equation with $160\times 2^7$ time steps, and for $X$, we use $160$ time steps. The processes $(Y, Z)$ are approximated iteratively, as per \eqref{eq:YZ}.
\newline\newline
\noindent \textbf{Problem specific settings:}\newline
For the deep multi-FBSDE method, we set $K=3$ and define $\psi_1(t,x,y,z)=A(C - x)$, $\psi_2(t,x,y,z)=-A(C - x)$, and $\psi_3(t,x,y,z)=0$. With this configuration, we incorporate both the doubling effect (through $\psi_2$) and the removal of the mean-reverting term (through $\psi_1$). Furthermore, since $\psi_3$ is identically zero, the original control problem is also explicitly included in the training.

The matrices used for the state equation are the same, up to $\sigma$, $R_x$ and $G$, as in \cite{andersson2023convergence, huang2025convergence}. In fact $R_x=\mathrm{diag}(5,1,5,1,5)$ is $G=\mathrm{diag}(1,5,1,5,1)$ in the implementation of \cite{andersson2023convergence}, but reads wrongly $R_x=\mathrm{diag}(25,1,25,1,25)$ and $G=\mathrm{diag}(1,25,1,25,1)$ in the presentation of that paper). Our setting reads
 \begin{align*}
    A&=\text{diag}([1,2,3,1,2,3]),\quad B=\begin{pmatrix}1&-1\\1&1\\0.5&1\\1&-1\\0&-1\\0&1\end{pmatrix},\quad C=\text{diag}([-0.2,-0.1,0,0,0.1,0.2]),\\
    \sigma&=\text{diag}([0.2,1,0.2,1,0.2,1]),\quad x_0=(0.1,0.1,0.1,0.1,0.1,0.1)^\top,\quad T=0.5,
\end{align*}
and
\begin{equation*}
        R_x=\text{diag}([25,1,25,1,25,1]),\quad R_u=\text{diag}([1,1]),\quad
    G=\text{diag}([1,25,1,25,1,25]).
\end{equation*}
In \cite{huang2025convergence}, the authors use the deep FBSDE method on their control problem, to approximate the solution of a BSDE derived from the stochastic maximum principle. However, this problem is structurally less complex because, at each time $t$, the control process in the associated BSDE is constant in the state space. It remains to understand to what extent this approach would generalize to more complex problems where the control process is non-constant.

\begin{figure}[htp]
\centering
\begin{tabular}{ccc}
          \includegraphics[width=80mm]{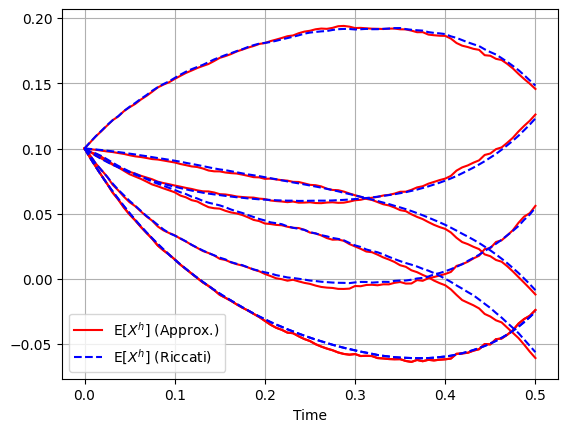}&  \includegraphics[width=80mm]{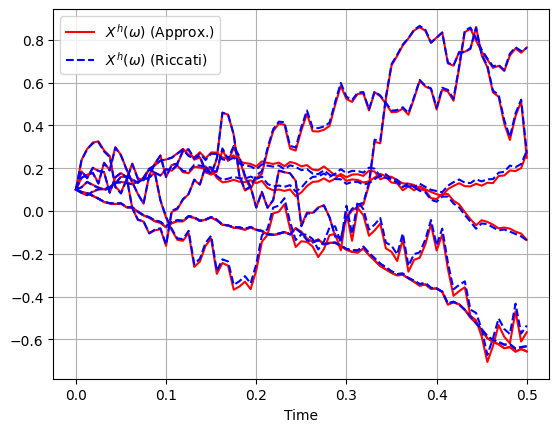}\\ \includegraphics[width=80mm]{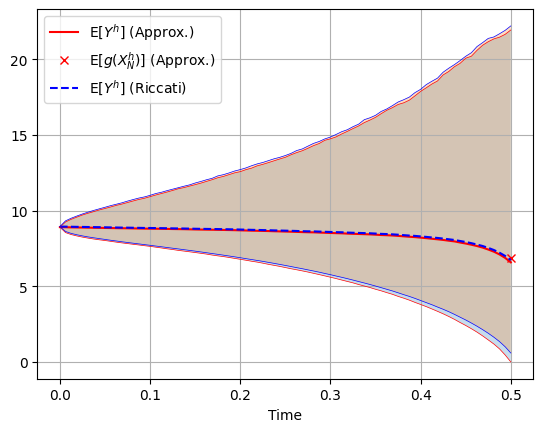}&
          \includegraphics[width=80mm]{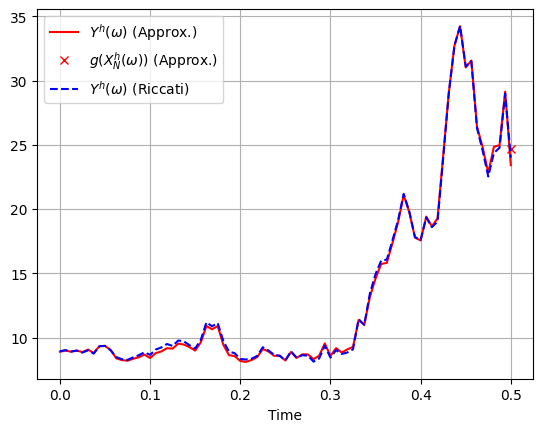}\\
          \includegraphics[width=80mm]{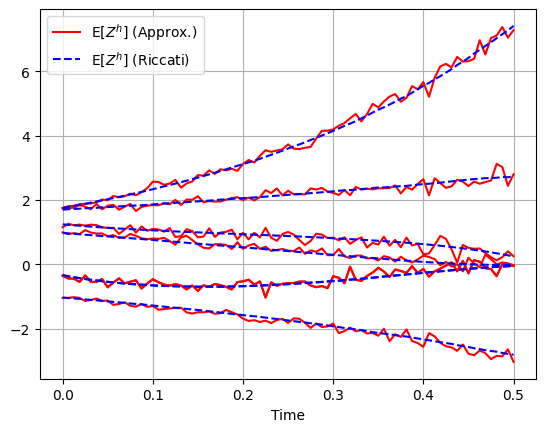}&   
         \includegraphics[width=80mm]{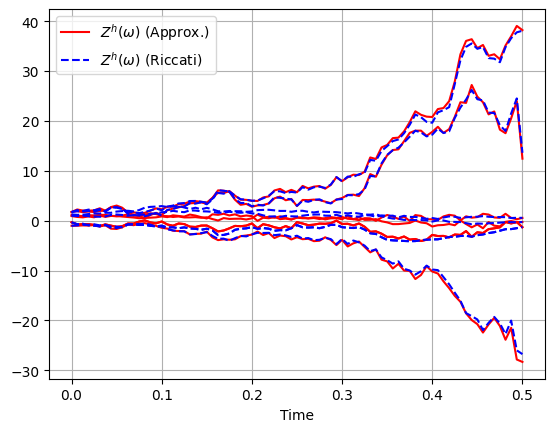}
\end{tabular}
\caption{Empirical mean of solutions (left) and a single solution trajectory (right) compared to the reference solution for the LQ-control problem from Section~\ref{lqgcp_num}. The shaded area represents an empirical credible interval for $Y$, defined as the area between the 5:th and the 95:th percentiles at each time point. We do not include credible intervals for $X$ and $Z$ in this figure to facilitate visualization. For $X$ and $Z$, we see one realization of each of the six components.}\label{fig:XYZ}
\end{figure}

\noindent \textbf{Results and discussion:}\newline
Figure~\ref{fig:XYZ} clearly demonstrates that our method yields accurate approximations for $(X,Y,Z)$, both pathwise and in distribution. We refrain from comparing these results with those obtained via the deep FBSDE method, as the latter produces completely incorrect approximations for this problem. The approximate initial condition obtained by the deep FBSDE method is around 28.44, whereas the true value is approximately 9 (see Figure~\ref{fig:MSE}). 

\begin{figure}[htp]
\centering
\begin{tabular}{ccc}
          \includegraphics[width=54mm]{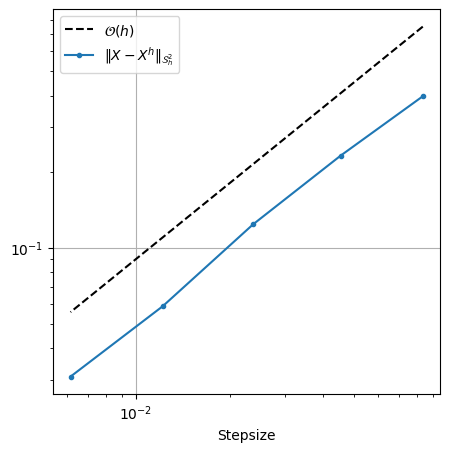}
        \includegraphics[width=54mm]{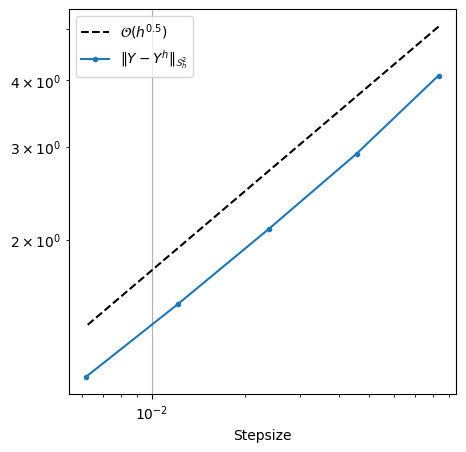}
          \includegraphics[width=54mm]{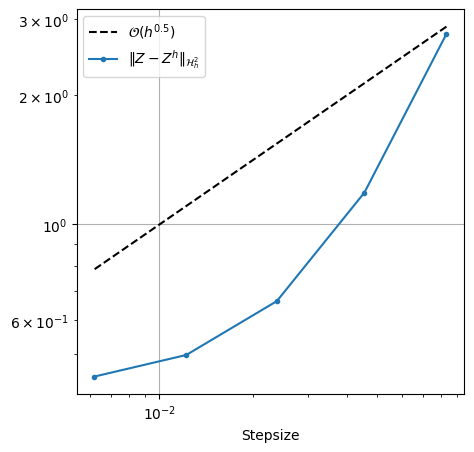}
          \end{tabular}
\caption{Empirical convergence plots for approximations of $(X,Y,Z)$ using the deep multi-FBSDE method phase II.}\label{fig:H2_XY_S2_Z}
\end{figure}

Part of the numerical experiments of this problem is the convergence to $Y_0$, shown in Figure~\ref{fig:L2_Y0}, as the time step $h$ tends to zero, or in the experiment for $N\in\{20,40,60,80\}$ time steps. To complete the empirical convergence study, Figure~\ref{fig:H2_XY_S2_Z} shows the convergence of $X$ and $Y$ with orders 1 and 0.5, respectively, and a decreasing error in $Z$ without a clear rate. We note that convergence orders 1 and 0.5 are the expected rates for forward SDEs with additive and multiplicative noise, respectively. The norms for the convergence are
\begin{equation*}
    \|A\|_{\mathcal{S}_{h,M}^2(\R^q)}=\max_{n\in\{0,1,\ldots,N\}}\bigg(\frac{1}{M}\sum_{m=1}^M\|A_n(m)\|^2\bigg)^{\frac12}, \quad \|A\|_{\mathcal{H}_{h,M}^2(\R^q)}=\frac{1}{N}\sum_{n=0}^{N-1}\bigg(\frac{1}{M}\sum_{m=1}^M\|A_n(m)\|^2\bigg)^{\frac12}.
\end{equation*}
Here, $A(m)=\{A_1(m),A_2(m),\cdot,A_N(m)\}$, $m=1,2,\ldots,M$, are \textit{i.i.d.} realizations of some adapted stochastic processes $A$ on the grid.

\subsection{FBSDEs from advection diffusion reaction PDEs}
In this section, we consider semilinear parabolic PDEs with quadratic gradient nonlinearity.

Let $\alpha,\beta\in\R^+$, $\gamma\in\R$, and $x_0\in\R^d$. For a general function $g\colon\R^d\to\R$ we define
\begin{equation}\label{PDE_parabolic}\begin{dcases}
\frac{\partial v}{\partial t} + \alpha\Delta v + 2\beta\|\mathrm{D}_x v\|^2 - \gamma v = 0, & (t,x) \in [0,T)\times \R^d,\\
v(T,x) = g(x), & x \in \R^d.
\end{dcases}
\end{equation}
The above PDE is structurally similar to a Hamilton--Jacobi--Bellman (HJB) equation. However, since $\beta>0$, it cannot be written in the form \eqref{PDE_CH}, does not have an associate stochastic control problem, and hence does not represent an HJB equation. We consider the equivalent FBSDE 
\begin{equation}\label{FBSDE_CHT_2}\begin{dcases}
X_t = x_0 + \int_0^t \beta Z_t \text{d} s + \sqrt{2\alpha} W_t,\\
Y_t = g(X_T) + \int_t^T\big(\beta\|Z_s\|^2-\gamma Y_s\big) \text{d} s - \sqrt{2\alpha} \int_t^T \langle Z_s, \text{d} W_s\rangle, \quad t\in[0,T].
\end{dcases}
\end{equation}
We again employ the Cole–Hopf transformation to derive a reference solution. Specifically, defining $u(t,x) = \log\big(\frac{\alpha}{2\beta} v(t,x)\big)$ leads directly to the simplified PDE
\begin{equation}\label{PDE_parabolic_CH}\begin{dcases}
\frac{\partial u}{\partial t}(t,x) + \alpha\Delta u - \gamma u\log{u} = 0, & (t,x) \in [0,T)\times \R^d,\\
u(T,x) = \exp\Bigg(\frac{2\beta}{\alpha}g(x)\Bigg), & x \in  \R^d.
\end{dcases}
\end{equation}
Back transformation yields $u(t,x)=\exp(\frac{2\beta}{\alpha}v(t,x))$. In the special case $\gamma=0$, we have for $t\in[0,T]$ \begin{equation}\label{MC_PDE_2}
    Y_t=-\frac{\alpha}{2\beta}\log{\bigg(E\Big[\text{e}^{-\frac{2\beta}{\alpha}g(X_t+\sigma\sqrt{T-t}\,\xi)}\big|X_t\Big]\bigg)},
\end{equation}
where $\xi$ is a standard $d$-dimensional normally distributed random vector.
\newline\newline
\noindent \textbf{Problem specific settings:}\newline
For the deep multi-FBSDE method, we opt for $K=2$, setting $\psi_1(t,x,y,z)=0$ and $\psi_2(t,x,y,z)=-\frac{\beta}{2}z$. We employ a terminal condition $g(x)= -|x_1-x_2|$, where $x=(x_1,x_2)^\top$. The discretization of the problem is set with $N=40$. In our first experiment, we set $\alpha=0.0315$, $\beta=0.6$, $x_0=(-0.1, 0.1)^\top$, and $\gamma=0$ and use \eqref{MC_PDE_2} to compute a reference solution. In our second experiment, we set $\alpha=0.0315$, $x_0=(-0.1, 0.1)^\top$, $\gamma=1$, and let $\beta$ vary between 0.00125 and 5.

\begin{figure}[htp]
\centering
\begin{tabular}{ccc}
\includegraphics[width=80mm]{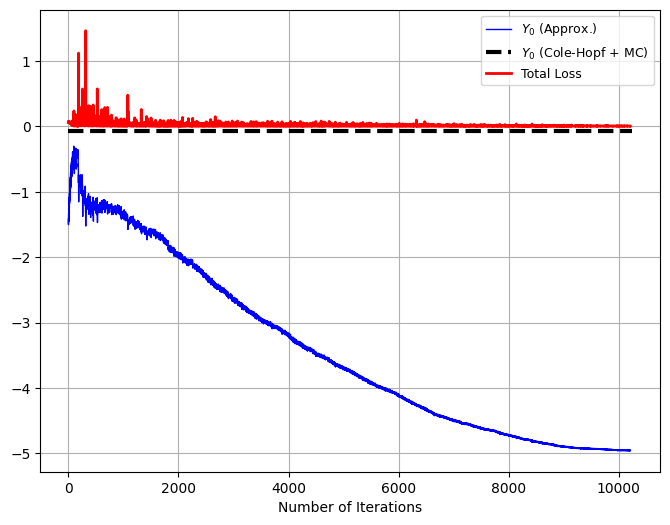}& \includegraphics[width=80mm]{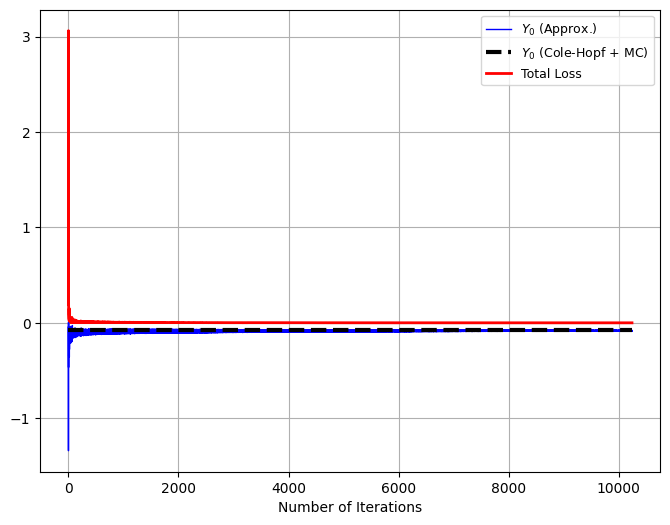}
\end{tabular}
\caption{Loss and approximation of $Y_0$ during training for the deep FBSDE method (left) and the deep multi-FBSDE method (right).
}\label{fig:FBSDE_CHT_2}
\end{figure}

\begin{figure}[htp]
\centering
\begin{tabular}{c}
          \includegraphics[width=160mm]{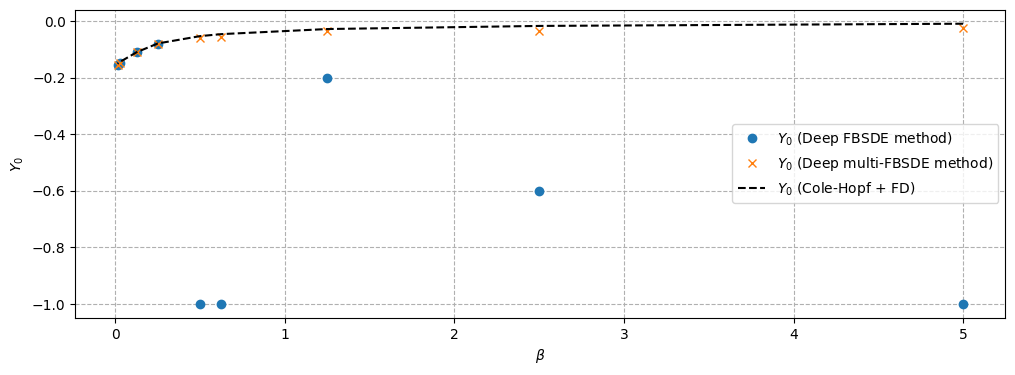}
          \end{tabular}
\caption{Approximate initial values for different values of $\beta$. Values of $Y_0$ smaller than $-1$ have been set to $-1$ for the sake of visibility.}\label{fig:Y0_b}
\end{figure}

For the first experiment, the reference solution is computed with a Monte Carlo approximation of \eqref{MC_PDE_2}. For the second experiment, we use a finite difference scheme to approximate \eqref{PDE_parabolic_CH}. More specifically, we approximate the solutions on a spatial domain $[-2,2]\times[-2,2]$ with artificial Dirichlet boundary conditions. 
\newline\newline
\noindent \textbf{Results and discussion:}\newline
Figure~\ref{fig:FBSDE_CHT_2} displays the approximate initial condition during training for both the deep FBSDE and deep multi-FBSDE methods.
The figure clearly illustrates that the deep FBSDE method yields inaccurate approximations, whereas the deep multi-FBSDE method produces highly accurate results. For the second experiment, we explore the accuracy of the deep FBSDE method and the deep multi-FBSDE method for different values of $\beta$. In Figure~\ref{fig:Y0_b}, we see that for $\beta>0.25$, the deep FBSDE fails to provide accurate approximations. The deep multi-FBSDE method provides accurate approximations for all $\beta$ tested. 

\section{Conclusions and potential future research directions}
In this paper, we introduced a two-phase deep-learning-based method, the deep multi-FBSDE method, for approximating coupled forward--backward stochastic differential equations (FBSDEs). Our approach focuses on scenarios where the standard deep FBSDE method \cite{han2018solving} does not converge. The main characteristic of our method is that its loss function is composed of the sum of multiple deep FBSDE loss functions, each originating from a FBSDE satisfying the same partial differential equation (PDE) as the original FBSDE, but with transformed drift and driver. The method was demonstrated on two stochastic control problems and two FBSDEs related to advection diffusion reaction PDEs, unrelated to stochastic control. For neither of the problems, the deep FBSDE method worked, but the deep multi-FBSDE method did. Overall, our results suggest that focusing on a family of equivalent FBSDEs offers a promising path to overcome convergence challenges in deep FBSDE methods, thus broadening the scope of problems in finance, control theory, and related areas where these methods can be applied effectively.

Although the method seems to work very well on challenging FBSDEs, we still do not understand the reason for this and why the deep FBSDE method fails. There is thus a need to both understand theoretically for what problems the deep FBSDE does not work and why, and why the deep multi-FBSDE method does work for these problems. From a practical perspective, the choice of $K$ and $\psi_1,\dots,\psi_K$ in our algorithm is currently somewhat \textit{ad hoc}. While this works for moderately complex problems, a more systematic approach to selecting $K$ and $\psi_1,\dots,\psi_K$ would be especially beneficial when dealing with high-dimensional systems or models where it is difficult to anticipate the impact of different choices. Developing such a framework constitutes a natural next step toward enhancing both the reliability and the applicability of the proposed method. Furthermore, extending the methodology to broader classes of PDE-driven problems and exploring theoretical guarantees for large-scale or high-dimensional settings would deepen our understanding of its robustness and versatility.

\section*{Acknowledgments}
Part of this work was carried out during a two-month postdoctoral visit of Kristoffer Andersson (KA) at the University of Utrecht, whose hospitality and support are gratefully acknowledged. KA also acknowledges financial support from RiBa 2022 Research Fund provided by the University of Verona. We very much thank Per Ljung att Saab for his careful reading and many comments that helped us improve the manuscript, to Karl Hammar at Saab and Chalmers for pointing out an error, and to Benjamin Svedung Wettervik at Saab for some interesting reflections.

\end{document}